\documentstyle[12pt]{article}
\newfont{\mib}{cmmib10 scaled 1200}
\newcommand{\dfrac}[2]{\frac{\displaystyle #1}{\displaystyle #2}}
\renewcommand{\phi}{\varphi}
\newtheorem{Lemma}{Lemma}[section]
\newtheorem{Theorem}{Theorem}[section]
\newtheorem{Proposition}{Proposition}[section]
\newtheorem{Corollary}[Proposition]{Corollary}
\newtheorem{Definition}{Definition}[section]
\newtheorem{Example}{Example}[section]
\begin{document}
\title{Bounds of automorphism groups of \\genus 2 fibrations}
\author{Zhi-Jie Chen\thanks{This work is realized under the
support of NSF grant \#DMS 9022140 while the author is visiting
Mathematical Sciences Research Institute. Also partly supported by
NSFC and K.~C.~Wong Education Foundation.}}
\date{(Second version)\\March 3, 1993}
\maketitle

It is well known that the automorphism group of a surface of general type is
 finite and bounded by a function of $K^2$ \cite{A}. Since then, several
authors worked on this subject and found better upper bounds of the group.
 Recently Xiao [11,12] has obtained a linear bound for this group.
 Hence it is natural to investigate the upper bounds for particular classes of
 surfaces. Here we are interested in the upper bounds of various
automorphism groups of
 surfaces with genus 2 pencils. As a first step, in the present paper, we will
 study the upper bounds of automorphism groups of genus 2 fibrations.

We always suppose that $S$ is a smooth projective surface over the complex
number field. A genus 2 fibration is a morphism $f:S\longrightarrow C$ where
$C$ is a projective curve such that a general fiber of $f$ is a smooth curve
of genus 2.

\begin{Definition}
An automorphism of the fibration $f:S\longrightarrow C$ is a pair of
automorphisms $(\tilde{\sigma},\sigma)$ where $\tilde{\sigma}\in\mbox{\rm Aut}
(S)$, $\sigma \in\mbox{\rm Aut}(C)$, such that the diagram
\begin{center}
\begin{picture}(40,26)(0,-26)
\multiput(5,-2.5)(0,-21){2}{\vector(1,0){30}}
\multiput(2.5,-5)(35,0){2}{\vector(0,-1){16}}
\multiput(1,-4)(35,0){2}{$S$}
\multiput(1,-25)(35,0){2}{$C$}
\put(1.5,-13){\makebox(0,0)[r]{$f$}}
\put(38.5,-13){\makebox(0,0)[l]{$f$}}
\put(20,-1.5){\makebox(0,0)[b]{$\tilde{\sigma}$}}
\put(20,-22.5){\makebox(0,0)[b]{$\sigma$}}
\end{picture}
\end{center}
commutes.
\end{Definition}

The automorphism group of fibration $f$ will be denoted by Aut($f$).
The main results of this paper are the following

\begin{Theorem}
Suppose $S$ is a surface of general type over the complex number field
with a relatively minimal genus 2 fibration $f:S\longrightarrow C$. Then
$$|\mbox{\rm Aut}(f)|\le 504K_S^2.$$
If $f$ is not locally trivial, then
$$|\mbox{\rm Aut}(f)|\le 288K^2_S.$$
In particular,
$$|\mbox{\rm Aut}(f)|\le\left\{\begin{array}{lll}
126K_S^2,&\qquad&\mbox{if }g(C)\ge 2;\\
144K_S^2,&&\mbox{if }g(C)=1;\\
120K_S^2+960,&&\mbox{if }g(C)=0.\end{array}\right.
$$
And these are all the best bounds.
\end{Theorem}

\begin{Theorem}
Suppose $S$ is a surface of general type over the complex number field
with a relatively minimal genus 2 fibration $f:S\longrightarrow C$,
$G$ is an abelian automorphism group of $f$. Then
$$|G|\le 12.5K^2_S+100.$$
This bound is the best.
\end{Theorem}

\begin{Theorem}
Suppose $S$ is a surface of general type over the complex number field
with a relatively minimal genus 2 fibration $f:S\longrightarrow C$,
$G$ is a cyclic automorphism group of $f$. Then
$$|G|\le\left\{\begin{array}{lll}
5K_S^2,&\qquad&\mbox{if $g(C)=1$, $K^2_S\ge12$;}\\
12.5K^2_S+90,&&\mbox{if $g(C)=0$.}
\end{array}\right.$$
These bounds are the best.
\end{Theorem}

\begin{Theorem}
Suppose $S$ is a minimal surface of general type over the complex
number field with a genus 2 fibration $f:S\longrightarrow C$ such that
$g(C)\ge2$, $G$ is a cyclic automorphism group of $f$. Then
$$|G|\le 5K^2_S+30.$$
\end{Theorem}

Theorem 0.1 will be obtained as a consequence of several propositions
in Section 3. In Section 4, we discuss abelian and cyclic automorphism
groups of the fibration $f$. The propositions proved there imply
Theorems 0.2, 0.3 and 0.4. We remark that Xiao \cite{X6} has obtained a
bound for abelian automorphism groups of $f$. Our theorem is an
improvement of his. Examples are given in Section \ref{FS} to show
that most of these bounds are the best possible.

\section{Preliminaries}

The surfaces with genus 2 pencils have been largely studied by many authors.
The facts we needed in this paper mostly appeared in [3, 6, 9, 10]. In
particular, Xiao's book \cite{X3} gave a systematic description of the
properties of genus 2 fibrations which are just what we needed here.
Unfortunately, this book has not been translated into English yet, hence it is
not available for most readers. For this reason, we will recall some materials
in this section.

Let $f:S\longrightarrow C$ be a relatively minimal fibration of genus 2,
$\omega_{S/C}=\omega_S\otimes f^*\omega_C^\vee$ the relative canonical sheaf
of $f$. For a sufficiently ample invertible sheaf $\cal L$ on $C$, the natural
morphism $f^*(f_*\omega_{S/C}\otimes{\cal L})\longrightarrow \omega_{S/C}
\otimes f^*{\cal L}$ defines a natural map $\Phi$:
\begin{center}
\begin{picture}(75,25)(0,-25)
\multiput(5,-2.5)(1.5,0){18}{\line(1,0){1}}
\put(32,-2.5){\vector(1,0){3}}
\put(2.5,-5){\vector(1,-1){15}}
\put(37.5,-5){\vector(-1,-1){15}}
\put(20,-22){\makebox(0,0)[t]{$C$}}
\put(2.5,-4){\makebox(0,0)[b]{$S$}}
\put(36,-4){\makebox(0,0)[lb]{$P=\mbox{\mib
P}\,(f_*\omega_{S/C}\otimes {\cal L})$}}
\put(20,-1.5){\makebox(0,0)[b]{$\Phi$}}
\put(7,-10.5){\makebox(0,0)[tr]{$f$}}
\put(32,-11.5){\makebox(0,0)[tl]{$\pi$}}
\end{picture}
\end{center}
$\Phi$ is called a relatively canonical map. By a succession of blow-ups, we
can obtain the following commutative diagram:
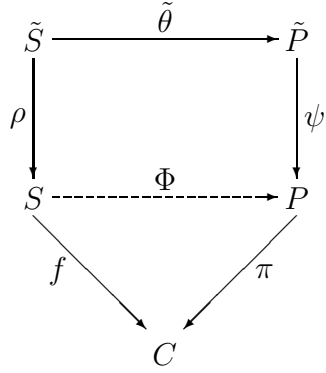
\begin{center}
\begin{picture}(40,49)(0,-25)
\multiput(5,-2.5)(1.5,0){18}{\line(1,0){1}}
\put(32,-2.5){\vector(1,0){3}}
\put(2.5,-5){\vector(1,-1){15}}
\put(37.5,-5){\vector(-1,-1){15}}
\put(20,-22){\makebox(0,0)[t]{$C$}}
\put(2.5,-4){\makebox(0,0)[b]{$S$}}
\put(37.5,-4){\makebox(0,0)[b]{$P$}}
\put(20,-1.5){\makebox(0,0)[b]{$\Phi$}}
\put(7,-10.5){\makebox(0,0)[tr]{$f$}}
\put(32,-11.5){\makebox(0,0)[tl]{$\pi$}}
\put(5,18.5){\vector(1,0){30}}
\multiput(2.5,16)(35,0){2}{\vector(0,-1){16}}
\put(2.5,17){\makebox(0,0)[b]{$\tilde{S}$}}
\put(37.5,17){\makebox(0,0)[b]{$\tilde{P}$}}
\put(1.5,8){\makebox(0,0)[r]{$\rho$}}
\put(38.5,8){\makebox(0,0)[l]{$\psi$}}
\put(20,19.5){\makebox(0,0)[b]{$\tilde{\theta}$}}
\end{picture}
\end{center}
where $\rho$ and $\psi$ are compositions of finitely many blow-ups,
$\tilde{\theta}$ is a double cover. Then we get the branch loci $\tilde R$ on
$\tilde P$ and $R$ on $P$ such that $\tilde R$ is the minimal even resolution
of $R$ (i.e. the canonical resolution of double cover). If $\cal L$ is
sufficiently ample, then all the singularities of $R$ must be located in one
of the 6 types of singular fibers defined by Horikawa \cite{H1}---0), I),
II), III), IV) or V).

$P$ is a relatively minimal ruled surface. We denote a section which
has the least self-intersection by $C_0$ such that $C_0^2=-e$. We will
use $F$ to denote both the fiber of $f$ or of $\pi$.

A singular point of the branch locus will be called {\it negligible} if
this point itself and all its infinitely near points are double points or
triple points with at least 2 different tangents. By minimal even resolution,
the inverse image of a negligible singular point is composed of $(-2)$-curves.
All other singular points are called {\it non-negligible}. The singular fiber
 of type 0) in the classification of Horikawa is nothing else but the fiber
which does not contain any non-negligible singular points.

The minimal even resolution $\psi:\tilde P\longrightarrow P$ can be decomposed
into $\tilde{\psi}:\tilde P\longrightarrow \hat P$ followed by $\hat{\psi}:
\hat{P}\longrightarrow P$, where $\tilde{\psi}$ and $\hat{\psi}$ are composed
respectively of negligible and non-negligible blow-ups. The image of $\tilde{R}
$ in $\hat{P}$ is denoted by $\hat{R}$.

If we take away all the isolated vertical $(-2)$-curves from the reduced
divisor $\hat{R}$, we get a new reduced divisor $\hat{R}_p$, which is called
the {\it principal part} of the branch locus $\hat{R}$. Then for any fiber $F$
of $\pi:P\longrightarrow C$, the second and third {\it singularity index} of
$F$, $s_2(F)$, $s_3(F)$, will be defined as follows:

If $R$ has no quadruple singularities on $F$, then $s_3(F)$ equals the number
of $(3\rightarrow3)$ type singularities of $R$ on $F$. Otherwise $s_3(F)$
equals the number of $(3\rightarrow3)$ type singularities of $R$ on $F$ plus
1. Hence $s_3(F)=0$ if and only if $R$ has no non-negligible singularities on
 $F$.

Let $\phi:\hat{R}_p\longrightarrow C$ be the natural projection induced by
$\pi\circ\hat{\psi}:\hat{P}\longrightarrow C$. Then the second singularity
index $s_2(F)$ of $F$ will be the ramification index of the divisor $\hat{R}_p$
on $f(F)$ with respect to the projection $\phi$. If $\hat{R}_p$ has
singularities (which must be negligible) on $F$, the singularity index $s_2(F)$
 can be calculated as follows.

For a smooth point $p\in\hat{R}_p\cap F$, the ramification index of $\phi$ at
$p$ can be defined as for an ordinary smooth curve. If $p\in\hat{R}_p\cap F$ is
 a singular point of $\hat{R}_p$, then the ramification index of $\phi$ at $p$
will be defined as the sum of ramification indices of the normalization of
$\hat{R}_p$ at the preimage of $p$ with respect to its projection to $C$ plus
the double of the influence to the arithmetic genus of $\hat{R}_p$ during its
normalization at the singular point $p$. If the normalization of $\hat{R}_p$
contains an isolated vertical component $E$, then the contribution of $E$ to
the ramification index of $\phi$ is equal to $2g(E)-2$.

As there are finite number of fibers $F$ with $s_i(F)\ne0$, we define the
$i$-th {\it singularity index} of $f$, $s_i(f)$, to be the sum of $s_i(F)$ for
all fibers, when $i=2$, 3. If we take away from the branch locus $R$ all the
fibers $F$ with odd $s_3(F)$, we obtain a divisor $R_p$ which is called the
{\it principal part} of $R$. Suppose that
$$R_p\sim -3K_{P/C}+nF,$$
where $K_{P/C}$ is the relative canonical divisor of $\pi$ and $\sim$
represents numerical equivalence. With these definitions, the formula
computing the relative invariants of a genus 2 fibration can be stated as
follows.

\begin{Theorem}[Xiao \cite{X3}]\label{ThmX}
Let $f:S\longrightarrow C$ be a relatively minimal fibration of genus 2. Then
$$K_{S/C}^2=K_S^2-8(g(C)-1)=\frac15s_2(f)+\frac75s_3(f)=2n-s_3(f),$$
$$\chi_f=\chi({\cal O}_S)-(g(C)-1)=\frac1{10}s_2(f)+\frac15s_3(f)=n-s_3(f).$$
\end{Theorem}

\section{Local cases}

We begin with a local fibration $f:S_{\Delta}\longrightarrow\Delta$ where $f$
is an analytic mapping onto the unit disk $\Delta$, $S_{\Delta}$ is a
2-dimensional analytic smooth manifold and the fibers of $f$ are projective
curves. We assume that the fiber of the zero is singular and all the fibers
over $\Delta^*=\Delta-\{0\}$ are smooth curves of genus 2.

Similarly, we have a commutative diagram
\begin{center}
\begin{picture}(40,49)(0,-25)
\multiput(5,-2.5)(1.5,0){18}{\line(1,0){1}}
\put(32,-2.5){\vector(1,0){3}}
\put(2.5,-5){\vector(1,-1){15}}
\put(37.5,-5){\vector(-1,-1){15}}
\put(20,-22){\makebox(0,0)[t]{$\Delta$}}
\put(2.5,-4){\makebox(0,0)[b]{$S_\Delta$}}
\put(38.5,-4){\makebox(0,0)[b]{$P_\Delta$}}
\put(20,-1.5){\makebox(0,0)[b]{$\Phi$}}
\put(7,-10.5){\makebox(0,0)[tr]{$f$}}
\put(32,-11.5){\makebox(0,0)[tl]{$\pi$}}
\put(5,18.5){\vector(1,0){30}}
\multiput(2.5,16)(35,0){2}{\vector(0,-1){16}}
\put(2.5,17){\makebox(0,0)[b]{$\tilde{S}_\Delta$}}
\put(37.5,17){\makebox(0,0)[b]{$\tilde{P}_\Delta$}}
\put(1.5,8){\makebox(0,0)[r]{$\rho$}}
\put(38.5,8){\makebox(0,0)[l]{$\psi$}}
\put(20,19.5){\makebox(0,0)[b]{$\tilde{\theta}$}}
\end{picture}
\end{center}
Denote the branch locus in $P_\Delta$ by $R_\Delta$. We also denote
the horizontal part of $R_\Delta$ by $R'_\Delta$, that is,
$$R'_\Delta=\left\{\begin{array}{lll}
R_\Delta-F_0,&\qquad&\mbox{if $R_\Delta$ contains $F_0$,}\\
R_\Delta&&\mbox{otherwise.}
\end{array}\right.$$
Let
$F_0=\pi^{-1}(0)$, $F_t=\pi^{-1}(t)$, $t\in\Delta^*$, and
$K_\Delta=\{\tilde{\sigma}\in\mbox{Aut}
(S_\Delta)|f\circ\tilde{\sigma}=f\}$. Any automorphism
$\tilde{\sigma}\in K_\Delta$ induces an automorphism $\sigma$ of
$P_\Delta$ satisfying $\pi\circ\sigma
=\pi$ and $\sigma(R_\Delta)=R_\Delta$. We denote the image of
$K_\Delta$ by $\bar{K}_\Delta\subseteq \mbox{Aut }P_\Delta$, then
$$|K_\Delta|=2|\bar{K}_\Delta|.$$

Note that any finite automorphism subgroup of {\mib P\,}$^1$ must be one of
the following:

\begin{center}
\begin{tabular}{llcl}
$G\subseteq \mbox{Aut({\mib P\,}$^1$)}$&&$|G|$&Number of points in an
orbit\\
Cyclic group&$Z_n$&$n$&1, $n$\\
Dihedral group&$D_{2n}$&$2n$&2, $n$, $2n$\\
Tetrahedral group&$T_{12}$&12&4, 6, 12\\
Octahedral group&$O_{24}$&24&6, 8, 12, 24\\
Icosahedral group&$I_{60}$&60&12, 20, 30, 60
\end{tabular}
\end{center}

For any $\sigma\in\bar{K}_\Delta$, its restriction to $F_t\cong
\mbox{\mib P\,}^1$, $\sigma|_{F_t}$, must preserve the set of 6 points
contained in $F_t\cap R_\Delta$. Hence $\bar{K}_\Delta$ can be
isomorphic to the following groups: $O_{24}$, $T_{12}$, $D_{12}$,
$D_6$, $Z_6$, $Z_5$, $D_4$, $Z_4$, $Z_3$, $Z_2$ and 1.

\begin{Lemma}\label{L1}
If $\bar{K}_\Delta\cong O_{24}$, $T_{12}$ or $D_{12}$, then $F_0$ is
contained  in $R_\Delta$ and
$R_\Delta$ has 6 ordinary double points on $F_0$.
 In this case, we have $s_2(F_0)=10$, $s_3(F_0)=0$.
\end{Lemma}

\noindent{\bf Proof.} Since $\bar{K}_\Delta\cong O_{24}$, $T_{12}$ or
$D_{12}$,
$R_\Delta\cap F_t$ ($t\in \Delta^*$) consists respectively of 6
vertices of a
regular octahedron, of 6 points corresponding to the centers of edges
of a regular tetrahedron, or of sixth roots of unit. These 6 horizontal
branches  of $R_\Delta$ cannot
intersect when $t\rightarrow 0$. But by assumption, $R_{\Delta}$ must
have some singularities, so $F_0$ is contained in $R_\Delta$.

Since $R_\Delta$ does not contain non-negligible singularities, one
has $s_3(F_0)=0$ and $R_\Delta=\hat{R}_\Delta=(\hat{R} _\Delta )_p$.
On $F_0$, $R_\Delta$ has 6 ordinary double points, the influence of
each double point to the arithmetic genus of $R_\Delta$ during its
normalization is equal to 1. The preimage of $F_0$ in the
normalization of $R_\Delta$ is a smooth vertical rational curve which
does not meet any other branches, so its contribution to the
index $s_2(F_0)$ is equal to $-2$. Therefore $s_2(F_0)=2\times
6+(-2)=10$.\hspace{\fill} $\Box$

\vspace{5mm}
We list the following useful lemmas, the proof is evident. Since local
equations are used for calculation of singularity indices, they are
given in simplified form, omitting some higher order terms. All the
non-negligible singularities here are canonical, i.e. defined by
Horikawa.

\begin{Lemma}\label{L2}
If $\bar{K}_\Delta\cong D_6$ and $R'_\Delta$ is not \'etale over
$\Delta$, then up to a coordinates transformation,

(1) The equation of $R'_\Delta$ is $(x^3-t^k)(t^kx^3-1)$, $k>0$. In
this case, $s_3(F_0)=0$ implies $s_2(F_0)\ge4$.

(2) The equation of $R'_\Delta$ is $(x^3-1)^2-t^k(x^3+1)^2$, $k>0$.
In  this case,
 we have $s_3(F_0)=0$, $s_2(F_0)\ge3$.
\end{Lemma}

\begin{Lemma}\label{L3}
If $\bar{K}_\Delta\cong Z_6$ and $R'_\Delta$ is not \'etale over
$\Delta$, then up to a coordinates transformation,
the equation of $R'_\Delta$ is $x^6-t^k$, $1\le k\le3$. If $k=3$, it
has a non-negligible singularity and $s_3(F_0)=1$, $s_2(F_0)=3$.
Otherwise $s_2(F_0)\ge5$.
\end{Lemma}

\begin{Lemma}\label{L4}
If $\bar{K}_\Delta\cong Z_5$ and $R'_\Delta$ is not \'etale over
$\Delta$, then up to a coordinates transformation,

(1) The equation of $R'_\Delta$ is $x(x^5-t^k)$, $k=1$, $2$. In this
case, $s_3(F_0)=0$, $s_2(F_0)\ge6$.

(2) The equation of $R'_\Delta$ is $x(t^kx^5-1)$, $k=1$, $2$. In this
case, $s_3(F_0)=0$, $s_2(F_0)\ge4$.
\end{Lemma}

\begin{Lemma}\label{L5}
If $\bar{K}_\Delta\cong D_4$ and $R'_\Delta$ is not \'etale over
$\Delta$, then up to a coordinates transformation,

(1) The equation of $R'_\Delta$ is $(x^2-1)((x-1)^2-t^k(x+1) ^2)
(t^k(x-1)^2-(x+1)^2)$, $k>0$. In this case, $s_3(F_0)=0$ implies
$s_2(F_0)\ge 6$.

(2) The equation of $R'_\Delta$ is $(x^2-1)(x^2-t^k)(t^kx^2-1)$,
$k>0$. In this case, we have $s_3(F_0)=0$, $s_2(F_0)\ge2$.
\end{Lemma}

\begin{Lemma}\label{L6}
If $\bar{K}_\Delta\cong Z_4$ and $R'_\Delta$ is not \'etale over
$\Delta$, then up to a coordinates transformation,
the equation of $R'_\Delta$ is $x(x^4-t^k)$, $k=1$, $2$. In this case,
we have $s_3(F_0)=0$, $s_2(F_0)\ge5$.
\end{Lemma}

\begin{Lemma}\label{L7}
If $\bar{K}_\Delta\cong Z_3$ and $R'_\Delta$ is not \'etale over
$\Delta$, then up to a coordinates transformation,

(1) The equation of $R'_\Delta$ is $(x^3-t^{k_1})(t^{k_2}x^3- a(t))$,
$k_1, k_2>0$, $a(0)\ne 0$. In this case,
$s_3(F_0)=0$ implies $s_2(F_0)\ge 4$.

(2) The equation of $R'_\Delta$ is $x^6+a(t)x^3+t^k$, $1\le k\le 3$.
 In this case, $s_3(F_0)=0$ implies $s_2(F_0)\ge5$.

(3) The equation of $R'_\Delta$ is $(x^3-b-t^{k_1})(x^3-b-
t^{k_2}a(t))$, $k_1, k_2>0$, $a(0)\ne0$ and $b\ne0$.
In this  case, we have $s_3(F_0)=0$, $s_2(F_0)\ge 6$.

(4) The equation of $R'_\Delta$ is $(x^3-t^k)(x^3-a(t))$, $1\le k\le3$,
$a(0)\ne0$. In
this case, we have $s_3(F_0)=0$, $s_2(F_0)\ge2$.

(5) The equation of $R'_\Delta$ is $((x-b)^2-t^ka(t))((x-b\omega)^2-
\omega^2t^k a(t))((x-b\omega^2)^2-\omega t^ka(t))$, $k>0$, $a(0)\ne0$,
$b\ne0$, $\omega=\exp(2\pi i/3)$. In this case, we have $s_3(F_0)=0$,
$s_2(F_0)\ge 3$.
\end{Lemma}

We summarize the results of Lemmas \ref{L2} through \ref{L7} in the
following table. Here we assume that $R'_\Delta$ has only negligible
singularities or ramifications on $F_0$.

$$\begin{array}{lccc}
\bar{K}_\Delta&|K_\Delta|&\multicolumn{1}{c}{s_2(F_0)}
&\multicolumn{1}{c}{|K_\Delta|/s_2(F_0)}\\
D_6&12&\ge3&\le4\\
Z_6&12&\ge5&\le2.4\\
Z_5&10&\ge4&\le2.5\\
D_4&8&\ge2&\le4\\
Z_4&8&\ge5&\le1.6\\
Z_3&6&\ge2&\le3\\
Z_2&4&\ge1&\le4\\
1&2&\ge1&\le2
\end{array}$$

\begin{Lemma}\label{L8}
If $R'_\Delta$ has only negligible singularities or ramifications on
$F_0$, then $|K_\Delta|/s_2(F_0)\le4$. Moreover, if $\bar
K_{\Delta}\cong Z_6$, $Z_5$, $Z_4$ or 1, then $|K_\Delta|/s_2(F_0)$
$\le2.5$.
\end{Lemma}

\section{Bounds of automorphism groups}

Let $G=\mbox{Aut}(f)$ be the automorphism group of the fibration of
genus two $f:S\longrightarrow C$. Then we have an exact sequence
$$\begin{array}{ccccccccc}
1&\longrightarrow &K&\longrightarrow &G&\longrightarrow
&H&\longrightarrow& 1,\\
&&&&(\tilde{\sigma},\sigma)&\mapsto&\sigma&
\end{array}
$$
where $H\subseteq \mbox{Aut}(C)$, $K=\{(\tilde{\sigma}, \mbox{id}) \in
G\}=\{\tilde{\sigma}\in\mbox{Aut}(S)|f\circ \tilde{\sigma} =f\}$. Thus
$$|G|=|K||H|.$$
The elements of $H$ are often regarded as transformations of the
fibers of $f$ or $\pi$.

\begin{Proposition}\label{P1}
If $f:S\longrightarrow C$ is a relatively minimal fibration of genus 2
with $g(C)\ge2$, then
$$|\mbox{\rm Aut}(f)|\le504K^2_S.$$
\end{Proposition}

\noindent{\bf Proof.} Since $|K|\le48$,
$|H|\le|\mbox{Aut}(C)|\le84(g(C)-1)$, we have
$$|G|=|K||H|\le4032(g(C)-1).$$

On the other hand, $K^2_{S/C}\ge0$ and the equality holds if and only
if $f$ is locally trivial. Hence
$$K_S^2\ge8(g-1)(g(C)-1)=8(g(C)-1),$$
and
$$|G|\le504K_S^2.$$

\vspace{-8mm}\hspace*{\fill}$\Box$

\begin{Proposition}\label{P2}
If $f:S\longrightarrow C$ is a relatively minimal fibration of genus 2
which is not locally trivial with $g(C)\ge2$, then
$$|\mbox{\rm Aut}(f)|\le126K_S^2.$$
\end{Proposition}

\noindent{\bf Proof.} Let $R'$ denote the horizontal part of the
branch locus $R$. If $R'$ is not \'etale over $C$, then by
Lemmas of Section 2, we have $|K|\le12$. As $|H|\le84 (g(C)-1)\le
10.5K_S^2$,
$$|G|\le12|H|\le126K_S^2.$$

Now assume that $R'$ is \'etale. Since $f$ is not
locally trivial, we must have $K_{S/C}^2>0$, i.e. either $s_3(f)>0$ or
$s_2(f)>0$. So $R$ must contain some fiber $F_0$. By Lemma~\ref{L1},
$s_3(F_0)=0$, $s_2(F_0)=10$. Let $p=f(F_0)$, $n=|H|$. Since $H$ is a
subgroup of Aut$(C)$, $H$ determines a finite morphism $\tau:C
\longrightarrow X=C/H$. Denote the ramification index of $p\in C$ with
respect to $\tau$ by $r$ and the other ramification indices by $r_i$.
Then the Hurwitz's theorem implies that
$$2g(C)-2=n(2g(X)-2)+n\sum\left(1-\frac1{r_i}\right).$$
As the $H$-orbit of the point $p$ has $n/r$ points, this implies that
$s_2(f)\ge10n/r$. Hence
\begin{eqnarray*}
K_S^2&\ge&\dfrac15s_2(f)+8(g(C)-1)=\dfrac{2n}r +4n\left[2g(X) -2+ \sum
\left(1-\dfrac1{r_i}\right)\right]\\
&=&4n\left[2g(X)-2+\dfrac1{2r}+\sum\left(1-\dfrac1{r_i}\right)\right].
\end{eqnarray*}

It is not difficult to see that the expression $2g(X)-2+1/2r+ \sum(1-
1/r_i)$ reaches its minimal value $2/21$ (under the condition
$2g(X)-2+ \sum(1-1/r_i)>0$) when $g(X)=0$, $r_1=2$, $r_2=3$, and
$r=r_3=7$. That is
$$K_S^2\ge\frac8{21}n=\frac8{21}|H|.$$
Thus
$$|G|\le48|H|\le126K_S^2.$$

\vspace{-8mm}
\hspace*{\fill}$\Box$

\vspace{5mm}
\noindent{\bf Remark.} It is not difficult to see that if $f$ is not
locally trivial with $g(C)\ge2$ and $|\mbox{Aut}(f)|=126K_S^2$, then
$|\mbox{Aut}(C)|=84(g(C)-1)$, $|\mbox{Aut}(F)|=48$ for any smooth fiber $F$ and
Aut$(f)\cong \mbox{Aut}(C)\times\mbox{Aut}(F)$. We will give an
example later. In this case, the fibration $f$ is of constant moduli
and $S/\mbox{Aut}(f)\cong{\mib F}_1$.

\begin{Lemma}\label{L31}
Let $S$ be a surface of general type which has a relatively minimal
genus 2 fibration $f:S\longrightarrow C$. If the third singularity
index $s_3(f)\ne0$, then
$$|\mbox{\rm Aut}(f)|\le \frac{60}7rK^2_{S/C},$$
where
$$r=\min_{\mbox{${\scriptstyle s_3(F)\ne0}$}}|\mbox
{\rm Stab}_Hf(F)|,$$
$\mbox{\rm Stab}_Hf(F)$ is the stabilizer of $f(F)$ in $H$.
\end{Lemma}

\noindent{\bf Proof.}
Let $F_0$ be a singular fiber such that $s_3(F_0)\ne 0$ and
$r=|\mbox{Stab}_H f(F_0)|$. Then
$$K^2_{S/C}\ge\frac75s_3(f)\ge\frac{7s_3(F_0)}{5r}|H|,$$
and we get
$$|G|=|K||H|\le\frac{r}{s_3(F_0)}\cdot\frac{60}7K^2_{S/C}
\le\frac{60}7rK^2_{K/C}.$$

\vspace{-8mm}\hspace*{\fill}$\Box$

\begin{Lemma}\label{L32}
Let $S$ be a surface of general type which has a relatively minimal
genus 2 fibration $f:S\longrightarrow C$. If the horizontal part $R'$
of the branch locus $R$ is not \'etale and has only negligible singularities or
ramifications, then
$$|\mbox{\rm Aut}(f)|\le20rK^2_{S/C},$$
where
$$r=\min_{\mbox{\scriptsize ${\scriptstyle F}$
 singular fiber}}|\mbox
{\rm Stab}_Hf(F)|.$$
\end{Lemma}

\noindent{\bf Proof.}
Let $F_0$ be a singular fiber with $r=|\mbox{Stab}_Hf(F_0)|$. Since
here
$$K^2_{S/C}\ge\frac15s_2(f)\ge\frac{s_2(F_0)}{5r}|H|,$$
we have
$$|G|=|K||H|\le\frac{r|K|}{s_2(F_0)}\cdot5K^2_{S/C}\le 20rK^2_{S/C},$$
by Lemma \ref{L8}.\hspace{\fill}$\Box$

\begin{Lemma}\label{L33}
Let $S$ be a surface of general type which has a relatively minimal
genus 2 fibration $f:S\longrightarrow C$. If the horizontal part $R'$
of the branch locus $R$ is \'etale, then
$$|\mbox{\rm Aut}(f)|\le24rK^2_{S/C},$$
where
$$r=\min_{\mbox{\scriptsize ${\scriptstyle F}$
 singular fiber}}|\mbox
{\rm Stab}_Hf(F)|.$$
\end{Lemma}

\noindent{\bf Proof.}
Let $F_0$ be a singular fiber with $r=|\mbox{Stab}_Hf(F_0)|$. By
assumption, we have $s_2(F_0)=10$. Hence
$$|G|=|K||H|\le\frac{r|K|}{s_2(F_0)}\cdot5K^2_{S/C}\le 24rK^2_{S/C}.$$

\vspace{-9.5mm}\hspace*{\fill}$\Box$

\vspace{5mm}
Let $\bar{K}$ denote the
subgroup in Aut$(P)$ which is induced by $K$. If $\sigma\in\bar{K}$,
then $\pi\circ\sigma=\pi$ and $\sigma(R)=R$. Let $K_1$ be a cyclic
subgroup of order $m$ of
$\bar{K}$, $Q=P/K_1$ be the quotient surface. Then $Q$ is a ruled
surface.  We have a commutative diagram
\begin{center}
\begin{picture}(40,25)(0,-25)
\put(5,-2.5){\vector(1,0){30}}
\put(2.5,-5){\vector(1,-1){15}}
\put(37.5,-5){\vector(-1,-1){15}}
\put(20,-22){\makebox(0,0)[t]{$C$}}
\put(2.5,-4){\makebox(0,0)[b]{$P$}}
\put(37.7,-3.5){\makebox(0,0)[b]{$Q$}}
\put(20,-1.5){\makebox(0,0)[b]{$\alpha$}}
\put(7,-11.5){\makebox(0,0)[tr]{$\pi$}}
\put(32,-11.5){\makebox(0,0)[tl]{$\pi'$}}
\end{picture}
\end{center}
Let $C_0$ and $C_\infty\sim C_0+eF$ be the reduced ramification
divisors of $K_1$. Let $C'_0$ be a section
of $\pi'$ with the least self-intersection ${C'_0}^2=-e'$, $F'$ be a general
fiber of $\pi'$. Then $\alpha^*C'_0=mC_0$,
$\alpha^*C'_\infty=mC_\infty$, $\alpha^* F'=F$ and $e'=me$.
Let $D=\alpha(R')$, $C'=C'_0+C'_\infty$ be the branch locus.
Then $C'\sim2C'_0+e'F'\sim -K_{Q/C}$.

\begin{Lemma}\label{L34}
Assume $\bar{K}\cong D_6$. If $R'$ is not \'etale and has only
negligible singularities or ramifications, then $f$ has more than one
$H$-orbit of singular fibers.
\end{Lemma}

\noindent{\bf Proof.}
Let $K_1$ be the unique cyclic subgroup of order 3 of $\bar{K}$. There
are 2 types of singular fibers as listed in Lemma \ref{L2}. Let $F_0$
be a singular fiber. Then the local equations of $D$ near $F_0$ are
(1) $(x-t^k)(t^kx-1)$, $k\le3$, (2) $(x-1)^2-t^k(x+1)^2$, $k>0$.
In case (1), $D$ meets $C'$ at 2 points in $F_0$. In
case (2), $D$ does not meet $C'$ in $F_0$.

If all the singular fibers of $f$ is of type (1), then $D$ is an
\'etale cover of $C$. This means that $a=e'$,
$C'\sim D$. Hence $DC'=0$ which is impossible because $D$ and
$C'$ meet in $F_0$.

If all the singular fibers of $f$ is of type (2), then $DC'=0$.
Hence $D\sim C'$ and $D(D+K_{Q/C})=0$. This means that $D$ is \'etale
over $C$. A contradiction.\hspace*{\fill}$\Box$

\begin{Lemma}\label{L35}
Assume $\bar{K}\cong D_4$. If $R'$ is not \'etale,
 then $f$ has more than one
$H$-orbit of singular fibers.

If $H$ is cyclic and $g(C)=0$, then
$$|\mbox{\rm Aut}(f)|\le 12.5K^2_{S/C}.$$
\end{Lemma}

\noindent{\bf Proof.}
In this case, there are 4 sections in $P$ which do not meet each
other. Hence $e=0$. $R'$ contains 2 of these sections and denoted by
$C_0$ and $C_\infty$. Let $K_1$ be a cyclic subgroup of $\bar K$ with
$C_0$ and $C_\infty$ as ramifications. Assume that there is only one
$H$-orbit of singular fibers. If these singular fibers are all of type
(1) in Lemma \ref{L5}, then the local equation of
$D=\alpha(R'-C_0-C_\infty)$ is $(x-t^k)(t^kx-1)$, namely, $D$ is
\'etale. Therefore $D\sim 2C'_0$, $DC'_0=DC'_\infty=0$, a contradiction.
If the singular fibers are of type (2) in Lemma \ref{L5}, then $D$ does
not meet $C'_0$ and $C'_\infty$. Hence $D\sim 2C'_0$, $D^2=0$, a
contradiction. That implies there are at least 2 $H$-orbits.

Now suppose $H$ is cyclic. Let $h=|H|$. We call an $H$-orbit {\it big}
if it
contains $h$ fibers. If there is a big $H$-orbit whose singular fibers
are of type (1), then $s_2(F_0)\ge6$, so $|G|\le (20/3)K^2_{S/C}$. If
$|G|>(20/3)K^2_{S/C}$, then the singular fibers in a big $H$-orbit
must be of type (2) with $k\le2$. Let $F_2$ and $F_3$ denote 2 fibers
fixed by $H$. Then at least one of them is of type (1). The structure
of types (1) and (2) implies that the normalization of
$D=\alpha(R'-C_0-C_\infty)$ is \'etale with respect to $\pi'$. Hence
$D$ must decompose into 2 isomorphic sections $D_1$ and $D_2$,
$D_1\sim D_2\sim C'_0+aF'$. Since both $D_1$ and $D_2$ meet $C'_0$ and
$C'_\infty$, $F_2$ and $F_3$ are all singular of type (1). Since
$D_1D_2=2a=kh$, $D_1C'_0=a=kh/2$. Hence the local equation of $R'$
near $F_2$ or $F_3$ is
$(x^2-1)((x-1)^2-t^{kh/2}(x+1)^2)(t^{kh/2}(x-1)^2-(x+1)^2)$. When
$h\ge6$, these are non-negligible singularities. If $F_i$ ($i=2$, 3)
is singular fiber of type I), $s_3(F_i)=2[(kh-2)/8]+1\ge(kh-1)/4$. If
$F_i$ is of type II), $s_3(F_i)=2[kh/8]\ge(kh-6)/4$. So
$$K^2_{S/C}\ge \frac15\times 2\times
h+\frac75\times\frac{h-6}4\times2=\frac{11}{10}h-\frac{21}5.$$
$$|G|=8h\le\frac{80}{11}(K^2_{S/C}+\frac{21}5)<12.5K^2_{S/C}.$$
If there are more than one big $H$-orbits, it can be similarly shown
that $|G|\le 12.5K^2_{S/C}$.\hspace*{\fill}$\Box$

\begin{Lemma}\label{L36}
Assume $\bar{K}\cong Z_3$. If $R'$ is not \'etale and has only
negligible singularities or ramifications and $f$ has only one
$H$-orbit of singular fibers, then
$$|\mbox{\rm Aut}(f)|\le 6rK^2_{S/C}.$$
where
$$r=\min_{\mbox{\scriptsize ${\scriptstyle F}$
 singular fiber}}|\mbox
{\rm Stab}_Hf(F)|.$$
\end{Lemma}

\noindent{\bf Proof.}
Let $K_1=\bar K$. If the singular fibers are of types (1) or (4) in
Lemma \ref{L7}, then $D\sim 2C'_0+aF'$ is \'etale. $D(K_{Q/C}+D)=0$
implies $a=e'$. Hence $D(C'_0+C'_\infty)=0$, a contradiction. If the
singular fiber $F_0$ is of type (5) with $k=1$, then $D$ is
irreducible and smooth near $F_0$. This implies $DC'_\infty\ne0$, a
contradiction. Therefore $s_2(F_0)\ge5$ for any singular fiber $F_0$.
So $|G|\le6rK^2_{S/C}$.\hspace{\fill}$\Box$

\begin{Lemma}\label{L37}
Assume $\bar{K}\cong Z_2$. If $R'$ is not \'etale
 and $f$ has only one
$H$-orbit of singular fibers, then
$$|\mbox{\rm Aut}(f)|\le 5rK^2_{S/C}.$$
where
$$r=\min_{\mbox{\scriptsize ${\scriptstyle F}$
 singular fiber}}|\mbox
{\rm Stab}_Hf(F)|.$$
\end{Lemma}

\noindent{\bf Proof.}
Let $F_0$ be a singular fiber. $|G|>5rK^2_{S/C}$ implies
$s_2(F_0)\le3$. We distinguish between 2 cases.

{\it Case I.\/} $R'$ contains $C_0$ and $C_\infty$. Then the local
equation of $R'$ near $F_0$ must be (1) $x(x^2-t)(x^2-a(t))$, $a(0)\ne
0$, $s_2(F_0)=3$; (2) $x((x^2-a^2)^2-t)$, $a\ne0$, $s_2(F_0)=2$. Let
$D=\alpha (R'-C_0-C_\infty)\sim 2C'_0+aF'$. If all the singular fibers
are of type (1), then $D$ is \'etale. This is impossible. If the
singular fibers are of type (2), then $D$ is irreducible and does not
meet $C'$. This is impossible.

{\it Case II.\/} $R'$ does not contain $C_0$ and $C_\infty$. Then the
local equation of $R'$ may be (1) $(x^2-t)(x^2-a(t))(x^2-b(t))$,
$a(0)b(0)\ne0$, $a(0)\ne b(0)$, $s_2(F_0)=1$; (2)
$(x^2-t)(ta(t)x^2-1)(x^2-b(t))$, $a(0)b(0)\ne0$, $s_2(F_0)=2$; (3)
$((x^2-a^2)^2-t)(x^2-b(t))$, $ab(0)\ne0$, $s_2(F_0)=2$; (4)
$((x^2-a^2)^2-t)(x^2 -tb(t))$, $b(0)\ne0$, $s_2(F_0)=3$. Let
$D=\alpha(R')\sim 3C'_0+aF'$. If $F_0$ is of type (1) or (2), then $D$
is \'etale and smooth. $D$ must be decomposed into 3 disjoint
components. This means $e'=0$, a contradiction. If $F_0$ is of type
(3) or (4), then $D$ is smooth. The ramification index
$D(D+K_{Q/C})=4a-6e'=|H|/r$. Hence $DC'=2a-3e'=|H|/2r$. This is a
contradiction because we have $DC'=0$ for type (3) and $DC'=|H|/r$ for
type (4). \hspace{\fill}$\Box$

\begin{Proposition}\label{P3}
If $S$ is a minimal surface of general type which has a genus 2
fibration $f:S\longrightarrow C$ with $g(C)=1$, then
$$|\mbox{\rm Aut}(f)|\le144K_S^2.$$
\end{Proposition}

\noindent{\bf Proof.} In this case, we have
$$K_S^2=K_{S/C}^2=\frac15s_2(f)+\frac75s_3(f)>0.$$
Thus either $s_3(f)>0$ or $s_2(f)>0$.

Let $j(C)$ be the $j$-invariant of the elliptic curve $C$. Let $m$
denote the number of points contained in a smallest $H$-orbit of $C$.
Since $H$ is a finite subgroup of Aut$(C)$, we have
$$m=\left\{\begin{array}{lll}
|H|/2&\qquad&\mbox{if $j(C)\ne0$, 1728,}\\
|H|/4&&\mbox{if $j(C)=1728$,}\\
|H|/6&&\mbox{if $j(C)=0$.}
\end{array}\right.$$

Since $r\le6$, by Lemma \ref{L31}, \ref{L32} and \ref{L33}, the
conclusion is immediate.\hspace{\fill}$\Box$

\begin{Proposition}\label{P4}
If $S$ is a surface of general type which has a relatively minimal
fibration of genus two $f:S\longrightarrow C$ with $g(C)=0$, then
$$|\mbox{\rm Aut}(f)|\le120(K_S^2+8).$$
Moreover, we have
$$|\mbox{\rm Aut}(f)|\le48(K_S^2+8)$$
for $K^2_S\ge33$, and when $K^2_S\le32$, there are only 4 exceptions.
\end{Proposition}

\noindent{\bf Proof.} In this case, we have
$$K_S^2+8=K_{S/C}^2=\frac15s_2(f)+\frac75s_3(f)>0.$$
Hence either $s_3(f)>0$ or $s_2(f)>0$.

{\it Case I.\/}
Assume that $R'$ is \'etale over $C$.
If $r\le5$, then by Lemma \ref{L33}
$$|G|\le24rK^2_{S/C}\le120(K^2_S+8).$$

If $r\ge6$, then $H$ must be cyclic or dihedral group. In this case,
there are at most 2 singular fibers. Hence $K^2_{S/C}\le4$ by Theorem
\ref{ThmX}. This means $S$ is not of general type[10, Theorem
4.2.5, p.90].

{\it Case II.\/} Assume that $R'$ is not \'etale. Then $f$ is a
fibration of variable moduli. Hence $f$ must contain more than 2
singular fibers (\cite{Beauville}). This implies $r\le5$. The
conclusion is obtained by Lemmas \ref{L31} and \ref{L32}.

In the preceding argument, we can see that $|G|\le48(K^2_S+8)$ holds
if $r\le2$. If $|G|>48(K^2_S+8)$, we must have $r>3$. Then $H$ is one
of $T_{12}$, $O_{24}$ or $I_{60}$.

 If $f$ has more than one $H$-orbit
of singular fibers, then
\begin{eqnarray*}
\frac{K^2_{S/C}}{|G|}&\ge&\frac1{5r}\left(\frac{s_2(F_0)}{|K|}+\frac{7s_3(F_0)}
{|K|}\right)+\frac1{5r_1}\left(\frac{s_2(F_1)}{|K|}+\frac{7s_3(F_1)}{|K|}
\right)\\
&\ge&\frac1{25}\times\frac14+\frac1{20}\times\frac14=\frac9{400}>\frac1{48}.
\end{eqnarray*}
Therefore $f$ has only one $H$-orbit.

If the singular fibers has non-negligible singularities, then by Lemma
\ref{L31}, $|G|\le (60/7)rK^2_{S/C}\le(300/7)K^2_{S/C} <48K^2_{S/C}$.
Suppose that the horizontal part
$R'$ of the branch locus has only negligible singularities or
ramifications, then by Lemmas \ref{L34}, \ref{L35}, \ref{L36} and
\ref{L37}, we have
$$|G|\le12.5rK^2_{S/C}.$$
Thus $|G|>48K^2_{S/C}$ implies that $r\ge4$ and $\bar K$ is $Z_6$
or $Z_5$. If $\bar K\cong Z_6$,
then $r=5$ and $H\cong I_{60}$. To ensure $|G|>48K_{S/C}^2$, we have
$s_2(F_0)=5$, i.e. $R=R'\sim -3K_{P/C}+nF$ is a smooth irreducible
divisor. As a multiple cover on $C$, the ramification index of $R$ is
equal to $R(R+K_{P/C})=12n$. On the other hand, this ramification
index is equal to $5\times(60/5)=60$, i.e. $n=5$. But
$2n=10=K^2_{S/C}\ne s_2(f)/5=12$, a contradiction.

If $\bar K\cong Z_5$, then $|G|>48K^2_{S/C}$ implies $s_2(F_0)=4$. In
this case $R=R'= C_0+R_1$ where $R_1\sim 5C_0+(n+3e)F$ is an smooth
irreducible divisor and $R_1C_0=0$, i.e. $n=2e$. Computing the
ramification index of $R_1$ we get $R_1(R_1+K_{P/C})=10n=4|H|/r$. This
implies $5r||H|$, a contradiction. Hence $|G|>48(K^2_{S}+8)$ implies
that $R'$ is \'etale over $C$. There are only finite number of
possibilities.
 We list the possible fibrations with $|G|>48(K^2_S+8)$ as
follows.
$$\begin{array}{cccccc}
H&r&|G|&K^2_S&|K|/(K^2_S+8)&|K|/K^2_S\\
I_{60}&5&2880&16&120&180\\
I_{60}&3&2880&32&72&90\\
O_{24}&4&1152&4&96&288\\
O_{24}&3&1152&8&72&144
\end{array}$$

In the Section \ref{FS} we will show their
existence.\hspace{\fill}$\Box$

\begin{Corollary}\label{Cor}
If $S$ is a minimal surface of general type which has a genus 2
fibration $f:S\longrightarrow C$ with $g(C)=0$, then
$$|\mbox{\rm Aut}(f)|\le288K^2_S.$$
\end{Corollary}

\noindent{\bf Proof.}
If $K^2_S\ge2$, then
$$48(K^2_S+8)<288K^2_S.$$
By Proposition \ref{P4} we need only check the 4 exceptional examples.
That fact leads
to the inequality $|G|\le288K^2_S$.\hspace{\fill}$\Box$

\section{Abelian automorphism groups}

Let $G\subseteq\mbox{Aut}(f)$ be an abelian group. Then it is well
known that $|K|\le12$.

\begin{Proposition}[{[7, Lemma 8]}]\label{AP1}
Let $f:S\longrightarrow C$ be a relatively minimal fibration of genus
2 with $g(C)\ge2$, $G$ is an abelian automorphism group of $S$, then
$$|G|\le6K^2_S+96.$$
\end{Proposition}

Let $\bar G\subseteq\mbox{Aut}(P)$ be the induced automorphism group
of a commutative group $G$, then
$$1\longrightarrow \bar K\longrightarrow \bar G\longrightarrow H
\longrightarrow 1.$$

\begin{Lemma}\label{Z3}
Assume that $\bar K\cong Z_3$, $g(C)=0$. Let $p\in C$ be a fixed point
of the cyclic group $H$, $F=\pi^{-1}(p)$. If there is a $\bar K|_F$-orbit
containing 3 points in $F$, then
$$s_2(F)\ge3|H|.$$
\end{Lemma}

\noindent{\bf Proof.} Since $p$ is a fixed point of $H$, the induced
action of $\bar G$ on $F$ forms a commutative subgroup $\bar
G|_F\subseteq \mbox{Aut}(F)\cong \mbox{Aut({\mib P}\,$^1$)}$. Since
$\bar G|_F$ stabilize this $\bar K|_F$-orbit, $\bar G|_F=\bar K|_F\cong
Z_3$, i.e. $H|_F=1$. Hence the local equation of $R'$ near $F$ has the
form $f(x^3,t^h)$ where $h=|H|$. Or explicitely, the local equation of
$R'$ are (3) $(x^3-b-t^{k_1h}a_1(t^h))(x^3-b-t^{k_2h}a_2(t^h))$; (5)
$((x-b)^2 -t^{kh}a(t^h))((x-b\omega)^2-\omega^2t^{kh}a(t^h))
((x-b\omega ^2)^2-\omega t^{kh}a(t^h))$, $b\ne0$. Thus
$s_2(F)\ge3h=3|H|$. \hspace{\fill}$\Box$

\begin{Proposition}\label{AP2}
If $S$ is a surface of general type which has a relatively minimal
fibration of genus two $f:S\longrightarrow C$ with $g(C)\le1$, $G$ is
an abelian automorphism group of $f$, then
$$|G|\le12.5(K^2_S+8).$$
\end{Proposition}

\noindent{\bf Proof.} It is well known that $H$ must be a cyclic group or
a dihedral group $D_4\cong Z_2\oplus Z_2$.

If $g(C)=1$ and that $H$ does not act
freely on $C$, then $|H|\le6$. Hence
$|G|\le72<12.5(K^2_S+8)$. If $g(C)=0$ and $H\cong D_4$, then
$|G|\le48$, the claim holds too. So we can assume that $H$ is a cyclic
group and that there exists a singular fiber $F_0$ with
$|\mbox{Stab}_Hf(F_0)|=1$.

{\it Case I.\/} Suppose that the horizontal part $R'$ of the branch
locus $R$ is \'etale over $C$. Then
$$|G|\le6K^2_{S/C}.$$

{\it Case II.\/} Suppose that $R'$ is not \'etale. If there is a big
$H$-orbit with $s_3(F_0)\ne0$,
then
$$K^2_{S/C}\ge\frac75s_3(f)\ge\frac75|H|,$$
so
$$|G|\le\frac{60}7K^2_{S/C}<12.5(K^2_S+8).$$

Now suppose that on the big $H$-orbits $R'$ has only negligible
singularities or ramifications. If $\bar K\cong Z_6$, $Z_5$, $Z_4$ or
1, then by Lemma \ref{L8}, we have
$$|G|\le\frac{|K|}{s_2(F_0)}\cdot 5K^2_{S/C}\le 12.5K^2_{S/C}\le
12.5(K^2_S+8).$$

Suppose that $\bar K\cong D_4$, $Z_3$ or $Z_2$ and that
$|G|>12.5(K^2_S+8)$. Then Lemmas \ref{L35}, \ref{L36} and \ref{L37}
implies that $f$ must have more than one $H$-orbit of singular fibers.
To ensure $|G|>12.5(K^2_S+8)$, $f$ cannot have more than one big
$H$-orbits. Thus we have $g(C)=0$. Lemma \ref{L35} excludes the case
of $\bar K\cong D_4$.

If $\bar K\cong Z_3$, then $s_2(F_0)\le2$. Hence $F_0$ must be of type
(4) of Lemma \ref{L7} with $k=1$. Taking
$K_1=\bar K$ we construct the quotient surface $Q=P/K_1$ as in \S3.
Then $D=\alpha(R')$ is \'etale near $F_0$. But $D$ cannot be \'etale.
Hence at least one of the $H$-stabilized fibers $F_2$ or $F_3$ is of
type (2) $k=1$ or type (5) $k=1$. Lemma \ref{Z3} excludes the case of
type (5). Suppose one of the $F_i$ is of type (2). Then $D\sim
2C'_0+aF'$ is irreducible and smooth. As a smooth double cover of
$C\cong\mbox{\mib P\/}^1$,
the ramification index of $D$ is at least 2. So $F_2$ and $F_3$ are
all of type (2). Then $DC'=D(D+K_{Q/C})=2(a-e')=2$, a contradiction.

If $\bar K\cong Z_2$, then $s_2(F_0)=1$. Hence the local equation of
$R'$ near $F_0$ is $(x^2-t)(x^2-a(t))(x^2-b(t))$, $a(0)b(0)\ne0$,
$a(0)\ne b(0)$. So $D=\alpha(R')$ is \'etale near $F_0$.

If $F_2$ and $F_3$ have no ramifications, the $D$ can be decomposed
into 3 components $D_i\sim C'_0+a_iF'$, $i=1$, 2, 3. These 3
components must meet each other on $F_2$ and $F_3$. So there exists at
least one point on $F_i$ where 3 components intersect. The local
eqution of $R'$ will be $(x^4+a(t)x^2+t^2)(x^2-t^2b(t))$. But as
$D_i(C'_\infty -C'_0)=e'$, we will have $|H|\le1$.

If $F_2$ or $F_3$ has ramifications, the equation of $R'$ near $F_i$
must be (1) $x^6-t$; (2) $(x^4-t)(t^ka(t)x^2-1)$, $a(0)\ne0$; (3)
$((x^2-a^2)^2-t)(x^2-t^k)$, $a\ne0$; (4) $((x^2-a^2)^2-t)(x^2-b(t))$,
$b(t)\ne0$. If $F_2$ is of type (1), then $D$ is irreducible and
smooth. As a smooth triple cover of $C\cong \mbox{\mib P\/}^1$, the
ramification index of $D$ is at least 4. Hence $F_3$ is of type (1) as
well. Let $D\sim 3C'_0+aF'$. Then $2DC'=D(D+K_{Q/C})=4$, impossible.
If $F_2$ is of type (2), then $D$ is smooth and cannot be irreducible.
$D$ has 2 components $D_1\sim2C'_0+aF'$ and $D_2\sim 2C'_0+bF'$. By
the same argument, we have $D_1C'=D_1(D_1+K_{Q/C})+2$. Hence
$D_1C'_0=0$ and $D_1D_2=0$. It is impossible.
\hspace*{\fill}$\Box$

\vspace{5mm}
Suppose that $G$ is a cyclic automorphism group of $f$. Similarly,
there is an exact sequence
$$1\longrightarrow K\stackrel{\alpha}{\longrightarrow}
G\stackrel{\beta}{\longrightarrow} H\longrightarrow 1$$
where $H\subseteq\mbox{Aut}(C)$, $K=\{(\tilde\sigma, \mbox{id})\in
G\}$. It is known that $|K|\le10$.

\begin{Lemma}\label{AG3-1}
Suppose that $f:S\longrightarrow C$ is a fibration and that $G$ is a
cyclic automorphism group of $f$. If there exists a point $p\in C$
such that

(1) $\sigma|_{f^{-1}(p)}\in K|_{f^{-1}(p)}$, for $\sigma\in G$ and
$\sigma$ stabilize $f^{-1}(p)$;

(2) $K\longrightarrow\mbox{\rm Aut}(f^{-1}(p))$ is injective.

Then $|K|$ and $|\mbox{\rm Stab}_H(p)|$ are coprime.
\end{Lemma}

\noindent{\bf Proof.} Let $H_1=\mbox{Stab}_Hp$, $F=f^{-1}(p)$. Let
$h=|H_1|$, $k=|K|$, $d=(h,k)$. Assume that $\sigma$ is a generator of
$\beta^{-1}(H_1)$. Then $\beta((\sigma^{k/d})^h)=1$ implies
$\sigma^{hk/d}\in K$. On the other hand, since $\sigma|_F\in K|_F$ by
(1), we obtain $(\sigma^{h/d})^k|_F=\mbox{id}_F$. Thus
$\sigma^{kh/d}=1$ by (2). This is impossible.\hspace{\fill}$\Box$

\begin{Proposition}\label{AP3}
If $S$ is a surface of general type which has a relatively minimal
fibration of genus two $f:S\longrightarrow C$ with $g(C)=1$, $G$ is
a cyclic automorphism group of $f$, then
$$|G|\le5K^2_S$$
for $K^2_S\ge 12$.
\end{Proposition}

\noindent{\bf Proof.} If $H$ does not act freely on $C$, then $|H|\le
6$. Hence $|G|\le60$ and the conclusion holds. Therefore we will assume
$H$ acts freely afterwards. So $G\cong K\times H$ and $G$ is cyclic if
and only if $(|K|,|H|)=1$.
We will discuss case by case.

{\it Case I.\/} Suppose that the horizontal part $R'$ of the branch
locus $R$ is \'etale over $C$. There exists a
singular fiber $F_0$ with $|\mbox{Stab}_Hf(F_0)|=1$. It is not
difficult to show that in this case
$$|G|\le5K_S^2.$$

{\it Case II.\/} Suppose that $R'$ is not \'etale.

({\it a\/}) $\bar K\cong Z_5$. Let $F_0$ be a singular fiber. The
local
equation of $R'$ near $F_0$ is (1) $x(x^5-t^k)$ or (2) $x(t^kx^5-1)$,
$k=1$, 2. We construct the quotient surface $Q=P/\bar K$ as in Section
3. $R'$ must contain one of the section $C_0$ or $C_\infty$.
We take this section away from $R'$, get a reduced divisor
$R_1$ with $R_1F=5$. Let $D=\alpha(R_1)$, then $D\sim C'_0+aF'$. Since
$DC'_0=0$, $a=e'=5e$. Thus $R_1\sim 5C_0+5eF$ and $R_1C_\infty=5e$.
Since the intersection number of $R_1$ and $F$ on the fiber $F_0$ is
equal to $k\le2$, the number of singular fibers must be a multiple of
5. But $|H|$ can not be divided by 5, hence
the singular fibers are located in different $H$-orbits. This means
$|G|\le5K^2_{S}$.

({\it b\/}) $\bar K\cong Z_4$. The local equation of $R'$ near a
singular fiber $F_0$ is $x(x^4-t^k)$, $k=1$, 2. We use the same
construction as in case ({\it a}). Then $R'$ must contain $C_0$ and
$C_\infty$. Let $R_1=R'-C_0-C_\infty$, $D=\alpha(R_1)$. Then $D\sim
C'_0+e'F'$. Similarly we deduce $R_1C_\infty=4e$. Since $|H|$ cannot
be even, there are more than one singular $H$-orbits. So
$|G|\le5K^2_{S}$.

({\it c\/}) $\bar K\cong Z_3$. If $f$ has only one $H$-orbit of
singular fibers and that $|G|>5K^2_S$, then $s_2(F_0)=5$, namely, the
local equation of $R'$ is $x^6+a(t)x^3+t$. Constructing the quotient
surface $Q=P/\bar K$, $D=\alpha(R')\sim 2C'_0+aF'$ is a smooth
irreducible curve and $r\ne|H|$. Since $DC'_0=0$, $DC'_\infty=|H|$, we
get $a=e'=3e=|H|$, i.e. $(|H|,|K|)=3$, a contradiction.

({\it d\/}) $\bar K\cong Z_2$. Lemma \ref{L37} ensures $|K|\le5K^2_S$.

({\it e\/}) $\bar K=1$. If $s_2(F_0)\ge2$, then
$|G|\le5K^2_{S/C}$. If $s_2(F_0)=1$, there is only one situation, i.e.
the local equation of $R'$ near $F_0$ is $(x^2-t)(x-a_1(t))(x-a_2(t))
(x-a_3(t))(x-a_4(t))(x-a_5(t))$, $a_i(0)\ne0$. Suppose that there is
only one singular $H$-orbit. Then $R'$ is a smooth 6-tuple
cover of $C$. The contribution of each singular fiber to the
ramification index equals 1. By Hurwitz formula,
$$2g(R')-2=6(2g(C)-2)+|H|.$$
So $|H|$ is even, a contradiction.
\hspace{\fill}$\Box$

\begin{Proposition}\label{genus0}
If $S$ is a surface of general type which has a relatively minimal
fibration of genus two $f:S\longrightarrow C$ with $g(C)=0$, $G$ is
a cyclic automorphism group of $f$, then
$$|G|\le12.5K^2_S+90.$$
\end{Proposition}

\noindent{\bf Proof.} If $R'$ is \'etale, we have $|G|\le5K^2_{S/C}$.
If there is a singular fiber in a big $H$-orbit with $s_3(F)>0$, then
$|G|\le(50/7)K^2_{S/C}$. Now assume that $R'$ has only negligible
singularities or ramifications in big $H$-orbits. If $\bar K\cong Z_4$
or 1, we have $|G|\le10K^2_{S/C}$ by Lemma \ref{L8}. When $\bar K\cong
Z_3$ or $Z_2$, if $f$ has only one $H$-orbit of singular fibers, then
Lemmas \ref{L36} and \ref{L37} ensure $|G|\le6K^2_{S/C}$. Otherwise,
by the proof of Proposition \ref{AP2}, $f$ has at least 2 big
$H$-orbits of singular fibers, hence $|G|\le10K^2_{S/C}$.

It remains the case of $\bar K\cong Z_5$. The proof of Proposition
\ref{AP3} tells us that if $f$ has only one big $H$-orbit of singular
fibers, then $f$ has another singular fiber which is stabilized by
$H$. By Lemma \ref{L4}, we have
$$K^2_{S/C}\ge \frac45(|H|+1),$$
so
$$|G|+10|H|\le12.5K^2_{S/C}-10=12.5K^2_S+90.$$

\vspace{-8mm}\hspace*{\fill}$\Box$

\vspace{5mm}
When $g(C)\ge 2$, we need the following lemma on the order of some
automorphisms of a curve. The proof of the lemma is just a
slight modification of that of the theorem of Wiman\cite{W}. For the
convenience of the reader, we include its proof here which is a
modified copy of the version given in [8, Lemma B].

\begin{Lemma}\label{AL1}
Let $H$ be a cyclic group of automorphisms of a curve $C$ of genus
$g\ge2$ such that the order of $|\mbox{Stab}_H(p)|$ is odd for any $p\in
C$.  Then
$$|H|\le3g+3.$$
\end{Lemma}

\noindent{\bf Proof.} Let $x$ be a non-zero element in $H$ with
maximal number of fixed points, $H'$ the subgroup of $H$ generated by
elements fixing all fixed points of $x$, $n$ the number of fixed
elements of $x$, $k$ the order of $H'$. Then $k$ must be odd.
Let $C'=C/H'$, $g'=g(C')$, and
let $\Sigma$ be the image of the set of fixed points of $H'$ on $C'$.
We have
\begin{equation}\label{eq}
2g-2=2kg'-2k+n(k-1).
\end{equation}
and the quotient group $H''=H/H'$ is a cyclic group of automorphisms
of $C'$ which satisfies the same condition imposed on $H$, i.e.
$|\mbox{Stab}_{H''} (p)|$ is odd for any $p\in C'$.

If $n=0$, then $g'\ge2$ and $|H|\le g-1$. If $n=2$, then because every
non-zero element of $H''$ induces a non-trivial translation on
$\Sigma$, we must have $|H''|\le2$, so $|H|\le2k$, then $|H|\le2g$ by
(\ref{eq})(note that $g'\ne0$ in this case). So we may assume $n\ge3$.

Suppose $g'=1$, $H''$ acts freely on $C'$. Considering the induced
action $H''$ on $\Sigma$, we see that $|H''|\le n$. So (\ref{eq})
gives $|H|\le 2g+n-2$. On the other hand, since $k\ge3$, (\ref{eq})
also gives $n\le g-1$, therefore
$$|H|\le3g-3$$
in this case.

Suppose $g'=1$, $H''$ does not act freely on $C'$, then $H''$ has a
fixed point. By assumption, $|H''|$ must be odd. This implies
$|H''|\le3$. So (\ref{eq}) gives
$$|H|\le2g+1.$$

Now suppose that $C'$ is a rational curve. Then the action of $H''$
has exactly 2 fixed points. So $|H''|$ must be odd.
If one of these two points is in $\Sigma$,
then $|H''|\le n-1$ in view of the action of $H''$ on $\Sigma$. Since
$|H''|$ is odd, we have $n\ge4$. So
$$|H|\le3g+3.$$

Suppose that $\Sigma$ and the two fixed points, $\xi$, $\eta$ of
$H''$ are disjoint. Let $H_1\subset H$ be the stabilizer of a point in
the inverse image of $\xi$. Then $[H:H_1]=k$. As the stabilizer of a
point in the inverse image of $\eta$ is also of index $k$ in $H$, we
see that any non-zero element in $H_1$ fixes exactly $2k$ points,
i.e., the inverse image of $\xi$ and $\eta$. Now we can replace $H'$
by $H_1$ and repeat the arguments above (note that the only conditions
we used are that non-trivial elements in $H'$ have same fixed point
set and that $H/H'$ acts faithfully on $\Sigma$). But then $\Sigma$
is composed of two orbits of $H''$, so $|H''|\le n/2$, whereby
$$|H|\le\frac32g+3$$
by (\ref{eq}).

At last we use induction on $g$. Suppose that $g'\ge2$ and $|H''|\le3g'+3$.
(\ref{eq}) gives
$$3g+3-(n-4)\frac{3(g-g')}{2g'-2+n}\ge|H|.$$
If $n\ge4$, we have done. If $n=3$, by assumption, we must have
$|H''|\le3$. Therefore
$$|H|\le\frac{3(2g+1)}{2g'+1}\le\frac35(2g+1) \le3g+3.$$

\vspace{-8mm}\hspace*{\fill}$\Box$

\begin{Proposition}\label{AP4}
If $f:S\longrightarrow C$ is a relatively minimal fibration of genus 2
with $g(C)\ge2$, $G$ is a cyclic automorphism group of $f$, then
$$|G|\le5K^2_S+30$$
for $K^2_S\ge48$.
\end{Proposition}

\noindent{\bf Proof.} (1) Assume that $|H|=4g(C)+2$ and $|K|=10$. Let
$g=g(C)$. By the theorem of Wiman (see the version given in [8, Lemma
B]), $C$ is a cyclic cover of {\mib P\,}$^1$ with ramification index
$r_1=2$, $r_2=2g+1$, $r_3=4g+2$ or $r_1=3$, $r_2=6$, $r_3=(4g+2)/3$.
In fact, these $r_i$ are the orders of Stab$_H(p)$ for $p\in C$. Since
$Z_{10}$ is a maximal cyclic automorphism subgroup of a smooth curve
of genus 2, by Lemma \ref{AG3-1} we have $(|\mbox{Stab}_H(p)|,|K|)=1$
if $f^{-1}(p)$ is a smooth fiber. But in case 1, $r_1$ and $r_3$ are
even, in case 2, $r_2$ and $r_3$ are even. So $f$ has at least
$(2g+10)/3$ singular fibers. By Lemma \ref{L4}, we have $s_2(F)\ge4$
for a singular fiber $F$. Hence
$$K^2_S-8(g-1)=K^2_{S/C}\ge\frac45
\cdot\frac{2g+10}3=\frac{8(g+5)}{15},$$
$$|G|=10|H|=40g+20\le \frac{75}{16}K^2_S+45\le5K^2_S+30$$
when $K^2_S\ge48$.

If $|K|\le8$ and $|K|$ is even, then by Lemma \ref{AL1} there exist
points $p\in C$ with $(|\mbox{Stab}_H(p)|,2)\ne1$. Hence
$K^2_S-8(g-1)=K^2_{S/C}\ge1$ and
$$|G|\le8|H|=32g+16\le 4K^2_S+44\le5K^2_S+30$$
when $K^2_S\ge14$.

If $|K|$ is odd, then $|K|\le5$. The inequality is immediate.

(2) Assume that $|H|$ is odd. By Lemma \ref{AL1}, we have
$$|H|\le3g+3.$$
So
$$|G|\le10|H|\le30g+30\le\frac{15}4K^2_S+60\le5K^2_S+30$$
when $K^2_S\ge24$.

(3) Assume that $|H|$ is even and $|H|<4g+2$. If $|K|=10$, $f$ must
have more than one singular fibers by Lemma \ref{L4}. So $K^2_S-8(g-1)
=K^2_{S/C}\ge2$. We get
$$|G|=10|H|\le40g\le5K^2_S+30.$$
If $|K|\le8$, it is not difficult to obtain this
inequality.\hspace*{\fill}$\Box$

\vspace{5mm}
It seems that this bound is not the best. In Section 5 we will give an
example to show there are infinitely many fibrations which has an
automorphism with order $3.75K^2_S+60$.

\section{Examples}\label{FS}

\begin{Example}\label{Ex1}
Fibration with $|G|=50K^2_S$.
\end{Example}

Let $C$ be a Hurwitz curve, i.e. $|\mbox{Aut}(C)|=84(g(C)-1)$, $F$ be
a curve of genus 2 with $|\mbox{Aut}(F)|=48$. Let $S=C\times F$,
$f=p_1:S\longrightarrow C$. Then $K^2_S=8(g(C)-1)$,
Aut$(f)\cong\mbox{Aut}(C)\times\mbox{Aut}(F)$,
$$|\mbox{Aut}(f)|=|\mbox{Aut}(C)|\cdot|\mbox{Aut}(F)|=504K^2_S.$$

\begin{Example}\label{Ex2}
Fibrations with $|G|=126K^2_S$ which is not locally trivial.
\end{Example}

Let $F=\mbox{\mib P\,}^1$. Let $p_1=0$, $p_2=\infty$, $p_3=1$,
$p_4=\sqrt{-1}$, $p_5=-1$, $p_6=-\sqrt{-1}$ be 6 points on $F$. Let
$C$ be a Hurwitz curve. Then $C$ has an $H$-orbit $\{ q_1,\dots,q_m\}$
which contains $m=12(g(C)-1)$ points. Let $P=C\times F$. Taking
$R=p^*_1(q_1+\dots+ q_m)+p_2^*(p_1+\dots+p_6)$ as the branch locus, we
construct a double cover of $P$. After desingularization, we get a
smooth surface $S$ with a genus 2 fibration $f:S\longrightarrow C$. By
computation, we obtain $K^2_S=32(g(C)-1)$, $|G|=48\times
84(g(C)-1)=126K^2_S$.

\begin{Example}\label{Ex3}
Fibrations with $|G|=144K^2_S$ and $g(C)=1$.
\end{Example}

Let $F$ and $p_1,\dots,p_6$ as Example \ref{Ex2}. Let $C$ be an
elliptic curve with $j$-invariant $j(C)=0$. Fix a $q_1\in C$, then the
order of the group of automorphisms Aut$(C,q_1)$ of $C$ leaving $q_1$
fixed is equal to 6. Let $H_1\cong Z_m\oplus Z_m$ be a subgroup of
translations of Aut$(C)$. Take an extension subgroup $H_1\subset
H\subset \mbox{Aut}(C)$ such that $H/H_1\cong \mbox{Aut}(C,q_1)$.
Then $|H|=6m^2$. Let $q_1,\dots,q_{m^2}$ be the orbit of $q_1$ under
$H$. Let $P=C\times F$. Using $R=p^*_1(q_1+\dots+q_{m^2})+p_2^*(p_1+
\dots+p_6)$ as the branch locus, we construct a double cover of $P$.
After desingularization, we get a smooth surface $S$ with a genus 2
fibration $f:S\longrightarrow C$. By computation, we get $K^2_S=2m^2$.
On the other hand, $|K|=48$ gives $|G|=288m^2=144K^2_S$.

\begin{Example}\label{Ex4}
Rational fibration with $|G|=120(K^2_S+8)$.
\end{Example}

Let $F$ and $p_1,\dots ,p_6$ as Example~\ref{Ex2}. Let $C=\mbox{\mib
P\,}^1$, $q_1,\dots,q_{12}$ be the 12 vertices of an icosahedron. Let
$P=C\times F$. Taking $R=p_1^*(q_1+\dots +q_{12})+p_2^*(p_1+\dots
+p_6)$ as the branch locus, we can construct a double cover of $P$.
After desingularization, we obtain a genus 2 fibration
$f:S\longrightarrow C$ with $K^2_S=16$, $|H|=60$, $|K|=48$,
$|G|=2880=120(K^2_S+8)$.

\begin{Example}\label{Ex5}
Rational fibrations with $|G|=48(K^2_S+8)$.
\end{Example}

Let $F$ and $p_1,\dots,p_6$ as Example~\ref{Ex2}. Let $C=\mbox{\mib
P\,}^1$, $q_1,\dots,q_m$ be the $m$-th roots of unit. Then use the
same construction as Example~\ref{Ex2}, we obtain a genus 2 fibration
with $K^2_S=2(m-4)$, $|K|=48$, $|H|=2m$, $|G|=96m=48(K^2_S+8)$.

\begin{Example}\label{Ex6}
Exceptional rational fibrations listed in the proof of Proposition
\ref{P4}.
\end{Example}

Using the same construction as Example~\ref{Ex2}, take
$q_1,\dots,q_{20}$ as the 20 vertices of a dodecahedron. We get a
fibration with $K^2_S=32$, $|G|=2880=90K^2_S$. If we take
$q_1,\dots,q_{6}$ as the 6 vertices of an octahedron, we get a
fibration with $K^2_S=4$, $|G|=1152=288K^2_S$. If we take
$q_1,\dots,q_{8}$ as the 8 vertices of a cube, we get a
fibration with $K^2_S=8$, $|G|=1152=144K^2_S$.

%\begin{Example}\label{Ex7}
%\end{Example}

\begin{Example}\label{Ex7}
Fibrations the order of whose abelian automorphism group is
$12.5(K^2_S+8)$.
\end{Example}

Let $x_0$, \dots, $x_{2m}$, $x_{2m+1}$ be the homogeneous coordinates
in $\mbox{\mib P\,}^{2m+1}$, $\mbox{\mib P\,}^{2m}$ be the hyperplane
defined  by $x_{2m+1}=0$. Let
$\phi:t \mapsto (1,t,\dots,t^{2m},0)$ be a $2m$-uple embedding of
$\mbox{\mib P\,}^1$ in $\mbox{\mib P\,}^{2m}$ and denote its image by
$Y$.  Then $Y$ is a rational
normal curve of degree $2m$. Let $X$ be the cone over $Y$ in
$\mbox{\mib P\,}^{2m+1}$ with vertex $P_0(0,0,\dots,0,1)$. Denote
$\eta=\exp(2\pi
i/10m)$. Then the automorphism $\sigma:(x_0,\dots,x_{2m+1})
\mapsto(x_0, x_1\eta,\dots,x_{2m}\eta^{2m},x_{2m+1})$ of
$\mbox{\mib P\,}^{2m+1}$ is of order $10m$. The automorphism
$\tau:(x_0,\dots,x_{2m+1})
\mapsto(x_0, \dots,x_{2m},x_{2m+1}\eta^{2m})$ of
$\mbox{\mib P\,}^{2m+1}$ is of order $5$. The cone $X$ is
stabilized  by these
automorphisms $\sigma$ and $\tau$. Take a hypersurface $H$ defined by
$x^5_0+x_{2m}^5+x_{2m+1}^5$ which is also stabilized by $\sigma$ and
$\tau$. Moreover,
$P_0\not\in H$. Now blowing up the cone $X$ at the vertex $P_0$ we
get a Hirzebruch surface $P=F_{2m}$ which has an automorphism
$\tilde\sigma$ of order $10m$ induced by $\sigma$ and an automorphism
$\tilde\tau$ of order 5 induced by $\tau$. The pull-back of
the intersection $H\cap X$ is a smooth divisor $R_1$ on $P$ which is
linearly equivalent to $5C_0+10mF$. Taking $R=R_1+C_0\equiv 6C_0+10mF$
which is a smooth even divisor and stabilized under $\tilde\sigma$ and
$\tilde\tau$, as
the branch locus, we can construct a double cover $S$ of $P$ which has
a natural genus 2 fibration $f:S\longrightarrow \mbox{\mib P\,}^1$.
Since $K_P\equiv -2C_0-(2m+2)F$,
$K_S^2=2(K_P+R/2)^2=8(m-1)$. The pull-back of $\tilde\sigma$ on $S$
can generate a cyclic automorphism subgroup $H$ of order $10m$. The
pull-back of $\tilde\tau$ on $S$ together with the hyperelliptic
involution of the fibration $f$ generates a cyclic automorphism
subgroup $K\cong Z_{10}$. As $H$ and $K$ commute, $G=KH\cong
Z_{10}\oplus Z_{10m}$ is an abelian automorphism group of $f$ with
order $|G|=100m=12.5(K^2_S+8)$.

\begin{Example}\label{Ex8+}
Rational fibrations which has an automorphism with order $12.5K^2_S+90$.
\end{Example}

Let $x_0$, \dots, $x_{2m}$, $x_{2m+1}$ be the homogeneous coordinates
in $\mbox{\mib P\,}^{2m+1}$, $\mbox{\mib P\,}^{2m}$ be the hyperplane
defined  by $x_{2m+1}=0$. Let
$\phi:t \mapsto (1,t,\dots,t^{2m},0)$ be a $2m$-uple embedding of
$\mbox{\mib P\,}^1$ in $\mbox{\mib P\,}^{2m}$ and denote its image by
$Y$.  Then $Y$ is a rational
normal curve of degree $2m$. Let $X$ be the cone over $Y$ in
$\mbox{\mib P\,}^{2m+1}$ with vertex $P_0(0,0,\dots,0,1)$. Denote
$\eta=\exp(2\pi
i/(50m-5))$. Then the automorphism $\sigma:(x_0,\dots,x_{2m+1})
\mapsto(x_0, x_1\eta^5,\dots,x_{2m}\eta^{10m},x_{2m+1}\eta)$ of
$\mbox{\mib P\,}^{2m+1}$ is of order $50m-5$. The cone $X$ is
stabilized  by this
automorphism $\sigma$. Take a hypersurface $H$ defined by
$x_0^4x_1+x^5_{2m}+x_{2m+1}^5$ which is also stabilized by $\sigma$
and $P_0\not\in H$. Now blowing up the cone $X$ at the vertex $P_0$ we
get a Hirzebruch surface $P=F_{2m}$ which has an automorphism
$\tilde\sigma$ of order $50m-5$ induced by $\sigma$. The pull-back of
the intersection $H\cap X$ is a smooth divisor $R_1$ on $P$ which is
linearly equivalent to $5C_0+10mF$. Taking $R=R_1+C_0\equiv 6C_0+10mF$
which is a smooth even divisor and stabilized under $\tilde\sigma$, as
the branch locus, we can construct a double cover $S$ of $P$ which has
a natural genus 2 fibration $f:S\longrightarrow \mbox{\mib P\,}^1$.
Since $K_P\equiv -2C_0-(2m+2)F$,
$K_S^2=2(K_P+R/2)^2=8(m-1)$. The pull-back of $\tilde\sigma$ on $S$
can generate a cyclic automorphism group $G_1$ of order $50m-5$. Since
$|G_1|$ is odd, $G_1$ and the hyperelliptic involution of the
fibration $f$ generate a
cyclic automorphism group $G$ of $S$. Therefore
$|G|=100m-10=12.5K_S^2+90$.

\begin{Example}\label{Ex8}
Fibrations which has an automorphism with order $5K^2_S$.
\end{Example}

Let $F=\mbox{\mib P\,}^1$. Let $p_1=0$, $p_k=\exp(2k\pi i/5)$,
$k=1,\dots,5$, be 6 points in $F$. Let $C$ be an elliptic curve,
$q_1,\dots,q_m$ be an orbit of a cyclic translation group $H\subseteq
\mbox{Aut}(C)$ of order $m$ where $m$ is an odd prime
different from 5. Then using the same construction as Example
\ref{Ex2}, we obtain a genus 2 fibration with $K^2_S=2m$, $K\cong
Z_{10}$. Let $G=K\times H\cong Z_{10m}$. $|G|=10m=5K^2_S$.

\begin{Example}\label{Ex9}
Fibrations which has an automorphism with order $3.75K^2_S+60$.
\end{Example}

Let $p_0=0$, $p_k=\exp(2k\pi i/3)$, $k=1$, 2, 3 be 4 points in
$C'=\mbox{\mib P\,}^1$. For any odd prime $m\ne 3$, 5, taking
$D=p_0+p_1+(m-1)p_2+(m-1)p_3$ as a branch locus, we can construct a
$m$-cyclic cover $\sigma:C\longrightarrow C'$. Then $g(C)=m-1$.
$H''=\{x\mapsto x\exp(2k\pi i/3)|k=1,2,3\}\cong Z_3$ is a cyclic
automorphism group of $C'$ which stabilizes the set
$\{p_0,p_1,p_2,p_3\}$. On the other hand, the Galois group $H'$ of the
$m$-cyclic cover $\sigma$ is isomorphic to $Z_m$. We can obtain an
extension
$$1\longrightarrow H'\longrightarrow H\longrightarrow
H''\longrightarrow 1$$
such that $Z_{3m}\cong H\subseteq\mbox{Aut}(C)$.

Let $q_0=0$, $q_k=\exp(2k\pi i/5)$, $k=1,\dots,5$, be 6 points in
$F\cong \mbox{\mib P\,}^1$. Let $P=C\times F$. Take
$R=p_2^*(q_0+q_1+\dots +q_5)$ as branch locus, we can construct a
double cover $\theta:S\longrightarrow P$ which is also a genus 2
fibration $f=p_1\circ\theta:S\longrightarrow C$. $F$ has a cyclic
automorphism group $K_1=\{y\mapsto y\exp(2k\pi i/5)|k=1,\dots,5\}\cong
Z_5$ which stabilizes the set $\{q_0,\dots,q_5\}$ and can be lift to
$P$. It is not difficult to see that we can get $K\cong Z_{10}$ by
adding the involution of the double cover. Then $G=K\times H\cong
Z_{30m}$ is a cyclic automorphism group of $f$ which satisfies
$$|G|=30m=30(g(C)+1)=\frac{15}4K^2_S+60,$$
because $K^2_S=8(g(C)-1)$.

\vspace{1cm}
Department of Mathematics

East China Normal University

Shanghai 200062

People's Republic of China

\end{document}